\documentclass[12pt]{amsart}
\usepackage{amscd}
\usepackage{amssymb}
\usepackage{amsmath}
\usepackage{amsthm}

\begin{document}
\newtheorem*{recall}{}
\newtheorem*{theorem}{ \bf Theorem}
\newtheorem{df}{ \sc Definition}[section]
\newtheorem{as}[df]{ \sc Assumption}
\newtheorem{ex}[df]{ \it Example}
\newtheorem{pr}[df]{ \sc Proposition}
\newtheorem{cl}[df]{ \sc Claim}
\newtheorem{thrm}[df]{ \sc Teorem\u a}
\newtheorem{cor}[df]{ \sc Corollary}
\newtheorem*{akn}{ \sc Aknowledgements}
\newtheorem{rem}{\sc Remark}
\newtheorem{pbl}[df]{ \sc Problem}
\newtheorem{cj}[df]{ \sc Conjecture}
\newtheorem{constr}[df]{ \sc Construction}
\newtheorem{lem}[df]{ \sc Lemma}
\newtheorem{sbs}[df]{\large\bf\hskip -6pt}
\newtheorem{eq}[df]{\rm}
\def\dim{{\rm dim}}
\def\deg{{\rm deg}}
\def\Tor{{\mathcal T}\!{\rm or}}
\def\codim{{\rm codim}}
\def\rank{{\rm rank}}
\def\det{{\rm det}}
\def\ann{{\rm Ann}}
\def\socle{{\rm socle\ }}
\def\hat{\widehat}
\def\bar{\overline}
\def\Pic{\rm Pic}
\def\Hilb{\rm Hilb}
\def\p{\phantom{\mpd{}^\mpd{}_\mpd{}}}
\def\mpr#1{\;\smash{\mathop{\hbox to 20pt{\rightarrowfill}}\limits^{#1}}\;}
\def\epi#1{\;\smash{\mathop{\hbox to 20pt{\rightarrowfill}\hskip
-15pt\rightarrow}\limits^{#1\,}}\;}
\def\mono{\lhook\joinrel\relbar\joinrel\rightarrow}
\def\mpd#1{\big\downarrow\rlap{$\vcenter{\hbox{$\scriptstyle#1$}}$}}
\def\mdp#1{\llap{$\vcenter{\hbox{$\scriptstyle#1$}}$}\big\downarrow}
\def\mup#1{\llap{$\vcenter{\hbox{$\scriptstyle#1$}}$}\big\uparrow}
\def\mpu#1{\big\uparrow\rlap{$\vcenter{\hbox{$\scriptstyle#1$}}$}}
\def\exseq#1#2#3{0\rightarrow #1 \rightarrow #2 \rightarrow #3 \rightarrow 0}
\def\Ac{{\mathcal A}}
\def\Bc{{\mathcal B}}
\def\Mc{{\mathcal M}}
\def\Hc{{\mathcal H}}
\def\Oc{{\mathcal O}}
\def\Ic{{\mathcal I}}
\def\Jc{{\mathcal J}}
\def\Tc{{\mathcal T}}
\def\C{{\mathbb C}}
\def\P{{\mathbb P}}
\def\Z{{\mathbb Z}}
\def\Q{{\mathbb Q}}
\def\A{{\mathbb A}}
\def\F{{\mathcal F}}
\def\Hom{{\mathcal H om}}
\def\banica{{B\u anic\u a}}

\def\Proof{\noindent\hskip4pt{\it Proof}.\ }
\def\qed{\hfill$\Box$\vskip10pt}

\title{Globally Generated Vector Bundles on ${\P}^3$ with $c_1=3$}
\author{Nicolae Manolache}
\date{}
\address{Institute of Mathematics, Romanian Academy, P.O.Box 1-764,
RO-014700 Bucharest, Romania}

\maketitle

\begin{abstract}
We classify the globally generated vector bundles on $\P^3$ whith the 
first Chern class $c_1=3$. The case $c_1=1$ is very easy, the case $c_1=2$ was done 
in \cite{SU}. The case of rank $2$ vector bundles was covered in \cite{Huh}.
The main tool used is Serre's theorem relating vector bundles of $\rank =2$ with 
codimension $2$ lci subschemes and its generalization for higher ranks, considered 
firstly by Vogelaar in \cite{Vog}.
\end{abstract}

\section{Introduction}
\label{intro}
We use freely the standard notation, namely that  introduced in \cite{FAC},
used (and extended  in \cite{EGA}) (cf. also \cite{HAG}).
Our base field is $\mathbb C$.

A globally generated nontrivial vector bundle $E$ on a a projective space has the 
first Chern class $c_1(E) > 0$. When $c_1(E)=1$, modulo a trivial summand, 
$E\cong \Oc (1)$ or $E \cong T(-1)$, $T$ being the tangent bundle. The case 
$c_1(E)=2$ is studied in \cite{SU}. We note that the case of rank $2$, globally 
generated vector bundle on $\P^n$ with $c_1=3$ is covered in \cite{Huh} and this study
was extended to $c_1 \le 5$ in \cite{CE}.

In this paper we address the case $c_1(E)=3$ on $\P ^3$. The classification of vector
bundles with $c_1=3$  on $\P^n$ will be considered in \cite{AM}. 

The aim of this paper is to prove:
\begin{theorem} On $\P^3$ the globally generated vector bundles with the first Chern 
class $c_1=3$, which are not direct sums of line bundles are as follows:

(i) None of rank $2$. From here the same result follows for $\P^n$, $n >3$. 

(ii) Those of rank $3$ are given by exact sequences of the shape:
\begin{align*}
\tag{1} 0 \to \Oc(-3) \to 4\Oc \to & E\to 0 \\
\tag{2} 0 \to \Oc(-2) \to 3\Oc\oplus \Oc(1) \to & E \to 0\\
\tag{3} 0 \to \Oc(-1) \to 2\Oc \oplus 2\Oc(1) \to & E \to 0 \\
\tag{4} 0 \to \Oc(-1) \to 3\Oc\oplus \Oc(2) \to & E \to 0\\
\tag{5} 0 \to \Oc(-2)\oplus \Oc(-1) \to 5\Oc \to & E \to 0\\
\tag{6} 0 \to 2\Oc(-1) \to 4\Oc\oplus \Oc(1) \to & E \to 0\\
\tag{7} 0 \to 3\Oc(-1) \to 6\Oc \to & E \to 0\\
\tag{8} 0 \to T(-2) \to 5\Oc\oplus \Oc(1) \to & E \to 0\ ,\\ 
\tag{9} 0 \to T(-2)\oplus \Oc(-1) \to 7\Oc \to & E \to 0
\end{align*}
where $T$ is the tangent bundle. 

(iii) Those of rank $ >3$  are extensions of those from (ii) by trivial bundles.
\end{theorem}

\begin{akn}
The idea of this paper originated in a discussion with Iustin Coanda about \cite{SU}, 
who pointed to me also some useful references. I want also to thank the Institute of 
Mathematics of the Oldenburg University for hospitality during the writing down of 
this paper, especially to Udo Vetter. Besides this I had no other support than the 
membership to the Institute of Mathematics of the Romanian Academy.
\end{akn}

\section{Preliminaries}
\label{pre}
We recall some general results, mostly  starting remarks also in 
\cite{SU}:\smallskip

\begin{itemize}
\item[$\mathbf A.$] (J.-P. Serre, cf. \cite{At}) If a vector bundle $E$ or rank 
$r > n$ on a smooth variety $X$ of dimension $n$ is globally generated, then $E$ 
contains a  trivial subbundle of rank $=r-n$ .

\item[$\mathbf B.$] If a  globally generated vector bundle $E$ or rank $r\le n$ 
on ${\P}^n$ has  $c_r=0$,  then $E$ contains a  trivial subbundle of rank $=1$ 
(cf. \cite{OSS}, chap. I, Lemma 4.3.2).

\item[$\bf C.$] If $E$ is a globally generated vector bundle on ${\P}^n$, then 
$E$ contains a trivial direct summand of rank $=\dim H^0(E^\vee )$ (cf. \cite{Ott}, 
Lemma 3.9).

\item[$\mathbf D.$] If $E$ is a vector bundle on $\P^n$, which is globally generated
and for which $H^0(E(-c_1))\not =0$, where $c_1$ is the first Chern class of $E$,
 then, modulo a trivial summand, $E\cong \Oc (c_1)$ (cf. \cite{S}, Prop.1).

\item[$\bf E.$] (Serre, Hartshorne's setting, cf. \cite{Ha}, Th. 1.1) There is a 
bijection  between the sets:

(i) triples $\{E,s,\varphi\}$, where $E$ is a rank 2 vector bundle on $\P^3$, $s \in 
H^0(E)$ is a section such that its zero-set $Y$ has codimension $2$, and 
$\varphi : \Lambda ^2E\cong L$.

(ii) pairs $\{Y,\xi \}$, where $Y$ is lci curve in $\P^3$ and $\xi :\omega _Y \cong
L(-4)|_Y =L\otimes \omega _{\P^3} \otimes \Oc _Y$.

\item[$\mathbf E'$.] In the above situation, $\Lambda ^2(E) \cong \Oc (c_1)$,
where $c_1$ is $c_1(E)$, $\deg\ Y =c_2(E)$, one has an exact sequence :
$$
0 \to \Oc \to E \to I_Y(c_1) \to 0
$$
and $\omega _Y\cong \Oc _Y(c_1-4)$. In particular: $2p_a(Y)-2=c_2(c_1-4)$

\item[$\mathbf E''.$] (cf. \cite{HM}, proof of 5.1.) If $E$ in $\mathbf E.$ is 
globally generated, then a general section of it has a smooth zero set.

\item[$\bf F.$] (Vogelaar (\cite{Vog}, here Arrondo's setting, cf. \cite{Arr}, 
Theorem 1.1) 

(i) If $E$ is a vector bundle of rank $r$ on $\P^n$ and $r-1$ sections 
$s_1,\ldots ,s_{r-1}$ of $E$  has as dependency locus a subscheme $Y$ of $\P^3$ 
which is lci, then one has an exact sequence:
$$
0 \to (r-1)\Oc _Y \xrightarrow{(s_1,\ldots ,s_{r-1})} E \to I_Y \otimes L \to 0 \ \ ,
$$
where $L=\Lambda ^rE$.
In this case, if $N$ is the normal sheaf of $Y$, then  
$\Lambda ^2 N \otimes L^\vee |_Y$ is generated by $r-1$ global sections.

(ii) Conversely, if $Y$ is a lci subscheme of codimension $2$ of $\P^3$ and $L$ is 
a line bundle on $\P^3$ with the property $H^2(L^\vee)=0$, if 
$\Lambda ^2 N \otimes L^\vee |_Y$ is generated by $r-1$ global sections, then there 
exists a vector bundle $E$ of rank $r$ such that $\Lambda ^rE=L$ and $E$ has $r-1$
global sections whose dependence locus is $Y$.

\item[$\mathbf G.$] If $\varphi : F \to G$ is a general morphism between vector 
bundles over a smooth variety and $\Hom (F,G)$ is globally generated, then the 
dependency locus $Y$ of $\varphi $ has codimension $\rank (G) - \rank (F) +1$ and 
is nonsingular outside a subset of codimension $ \ge \rank (G) - \rank (F) +3$ 
(cf. \cite{B}, 4.1 or \cite{Ch} for a more general result and \cite{PS}, 
\cite{Sa} for the case of arithmetically Cohen-Macaulay curves in $\P^3$).

\end{itemize}
\qed

According to {\bf A}, any vector bundle on $\P^n$ of rank $ >n$ is an extension of 
a rank $n$ vector bundle by a trivial one, so modulo extensions by trivial bundles,
it is enough to find the globally generated vector bundles of rank at most $3$
with $c_1=3$.

\section{Rank 2 Vector Bundles}
\label{Srank2}

This case is covered in \cite{Huh}. We give nevertheless other arguments.

If $E$ is a globally generated vector bundle of rank 2  on $\P^3$, nonsplittable, 
with $c_1(E)=3$, we may assume $H^0(E(-3))=0$, according to 
$\mathbf D$. A general section of it has a smooth zero set $Y$ which according to 
$\mathbf E'$ is a union of conics. 

The exact sequence:
\begin{equation}\label{rank2}
  0 \to \Oc \to E \to I_Y(3) \to 0 \hfill
\end{equation}
shows that the ideal of $Y$ should be generated by cubics and $\deg (Y) < 9$.

Were $Y$ a conic, $E$ would be decomposable: 
$E\cong \Oc (1) \oplus \Oc (2)$. This corresponds to $H^0(E(-2))\not = 0$.

I. {\bf Assume} $H^0(E(-2))=0$, $H^0(E(-1)) \not =0$. Denote $F:=E(-2))$. One has 
$c_1(F)=-1$ and $H^0(F)=0$, so F is a normalized stable vector bundle. Denote 
$c_2(F):=2n$ (recall that $c_1c_2 = 0 \mod 2$) and so $c_2(E)=2n+2$. As 
$H^0(I_Y(2)) \not =0$, $Y$ is contained in a complete intersection of type $(2,3)$
and so has degree at most 6, i.e. $n\le 2$.

(a) $n=1$. In this case $F$ is stable with $c_1=-1$, $c_2=2$. The moduli space of 
such vector bundles was studied in \cite{HaS} and \cite{Ma1}, where it was shown
that one has an exact sequence:
$$
0 \to \Oc \to F(1) \to I_Z(1) \to 0 \  ,
$$
with $Z$ a double structure on a line. The homogeneous ideal of $Z$ is 
$(ax+by,x^2,xy,y^2)$, where the ideal of the line $L$ which is the support of $Z$
is $(x,y)$, and $a$, $b$ are homogeneous polynomials in the projective coordinates
on $L$, of degree $3$, without common zeros. As $I_Z(2)$ is not globally generated,
E is neither.

(b) $n=2$. In this case $F$ is stable with $c_1=-1$, $c_2=4$. The moduli space of 
such vector bundles was studied in \cite{BM}, where it was shown that for those $F$ for 
which $H^0(F(1)) \not =0$ one has an  exact sequence:
$$
0 \to \Oc (1) \to F(2)\ (\cong E) \to I_Z(2) \to 0 \  ,
$$
where $Z$  is a lci curve of degree $4$, with $\omega _Z \cong \Oc _Z(-3)$
and $Z$ is either a double structure on a conic, with the ideal of the shape
$(ax+bq,x^2,xq,q^2)$, (if the ideal of the conic is $(x,q)$), with $\deg(a)=4$,
$\deg(b) =3$, or a union of to double lines of the shape from (a) above. 
In both cases $E$ is not globally generated, because $I_Z(2)$ is not.

II. {\bf Assume} $H^0(E(-1))=0$. Then  in the exact sequence \eqref{rank2}, $Y$ must
be a disjoint union of conics, not contained in a quadric, with the ideal generated 
by cubics. We show in the following that this situation cannot occur.
Indeed, $Y$ should be a union of at least $3$ disjoint conics. We show that  
$Y$ cannot lie in a complete intersection of two cubics, where at least one of them 
can be assumed to be smooth.
If $Y$ would lie in $Z=$ complete intersection of type $(3,3)$, let $Y'$ the liaison 
residue of $Y$ in $Z$. Then one has the exact sequence:
$$
0 \to \omega _Y \otimes \omega^{-1} _Z\ (=\Oc _Y(-3)) \to \Oc _Z \to \Oc _ {Y'} \to 0
$$
It follows that the Cohen-Macaulay curve $Y'$ has Hilbert polynomial $3t+12$, when
$Y$ consists of $3$ disjoint conics, and has Hilbert polynomial $t+19$ when $Y$
consists of $4$ disjoint conics. The second case is clearly absurd while the only 
curves locally Cohen-Macaulay  of degree 1 are the lines.
(A general result is known: the pairs $(d,p_a)$ of degree $d$ and arithmetical genus
$p_a$ for which the Hilbert scheme of locally Cohen-Macaulay curves $\Hilb (d,p_a)$ 
is not void were determined by several authors (cf. \cite{HaCM}. \cite{OS}, 
\cite{Sa}), namely: $d > 0, p_a = \frac{1}{2}(d-1)(d-2)$ or $d > 1, p_a \le 
\frac{1}{2} (d-2)(d-3)$.) Curves with Hilbert polynomial $3t+12$ do exist. 
Such a curve must have at most one non-reduced component. It can be only a line.
So $Y'$ contains a double line or is a triple line.
In our case we have the supplementary condition that $Y'$ is contained in a complete 
intersection of two cubic surfaces, one of them smooth. We show in the next section 
that on a smooth cubic surface there are double lines only of Hilbert polynomial 
$2t+3$ and triple lines only of Hilbert polynomial $3t+6$ (see in the next section 
the exact sequences ($\delta$) and ($\tau$)). So $Y'$ cannot have Hilbert polynomial
$3t+12$.

The assertion about $\P^n$, $n >3$ follows from the splitting criterium of Horrocks 
(cf. \cite{Ho} or \cite{OSS})

\section{Rank $\ge  3$ Vector Bundles}
\label{Srank3}

We use {\bf F.} Observe $\Lambda _Y:=\Lambda ^2N\otimes L^\vee |Y \cong \omega _Y 
\otimes \omega^{-1} _{\P^3} \otimes L^ \vee |Y$. We shall not restrict the analysis 
to $\rank = 3$, but we shall emphasize the quotients of degree $3$. 
By {\bf D.}, if $E$ is a nonsplittable globally generated vector bundle of rank $r$ 
with $c_1=3$, we may assume $H^0(E(-3))=0$.

If $Y$ is the dependency locus of general $r-1$ sections of $E$ then $Y$ is smooth
and one has:
\begin{equation}\label{rank3}
  0 \to (r-1)\Oc \to E \to I_Y(3) \to 0 \hfill
\end{equation}

In this case $\Lambda _Y \cong \omega _Y (1)$

(A) {\em Assume $H^0(E(-2))\not =0$.} Then $Y$ in \eqref{rank3} is a plane curve.
Assume $\deg (Y)=d$. Then $\Lambda _Y =\Oc _Y(d-2)$. If $\Lambda _Y$ is generated
by $r-1$ global sections (we do not pay attention in taking r minimal), then 
the minimal resolution of $I_Y$ gives  a resolution of $E$:
$$
0 \to \Oc (2-d) \to (r-1)\Oc \oplus \Oc (3-d) \oplus \Oc (2) \to E \to 0
$$
Here $E$ is globally generated for $d\le 3$. For $d=1$ one gets a contradiction,
for $d=2$ one gets $E \cong (r-2)\Oc \oplus \Oc (1) \oplus \Oc (2)$ and for $d=3$
one gets:
$$
0 \to \Oc (-1) \to r\Oc \oplus \Oc(2) \to E \to 0 \ \ \ .
$$
For $r \ge 4$, $E$ is an extension by a trivial bundle of a vector bundle given by
an exact sequence:
$$
0 \to \Oc (-1) \to 3\Oc \oplus \Oc (2) \to E' \to 0
$$

(B) {\em Assume $H^0(E(-2)) =0$, $H^0(E(-1)) \not =0$.} Then $Y$ is contained in a 
complete intersection of type $(2,3) $, and is not contained in a plane.

{\bf 1.} {\bf  Assume $Y$ in \eqref{rank3} connected.} We have to analyze the cases 
$\deg (Y) =3,4,5,6$

(i) When $\deg (Y)=3$, $Y$ is a rational cubic  and it is easy to see that $E$ 
should be, modulo trivial summands, $3\Oc (1)$.

(ii) When $\deg (Y)=4$, $Y$ is either a rational quartic, or a complete intersection
of type $(2,2)$. The last case gives that $E$, modulo trivial summands, is defined 
by an exact sequence:
$$
0 \to \Oc (-1) \to 2\Oc \oplus 2\Oc (1) \to E \to 0
$$
The shape of a minimal resolution of a rational quartic in $\P ^3$ is:
$$
0\to \Oc (-5) \to 4 \Oc (-4) \to 3\Oc (-3)\oplus \Oc (-2) \to I_Y \to 0
$$
This, together with \eqref{rank2}, gives the resolution of $E$:
$$
0\to \Oc (-2) \to 4\Oc (-1) \to (r+2)\Oc \oplus \Oc (1) \to E \to 0
$$
$Y$ being an embedding of $\P^1$ in $\P^3$, $\Lambda _Y$ pulled back to $\P^1$ is
$\Oc _{\P^1}(2)$, whence is generated by $2$ global sections. This shows that 
$r \ge 3$.
The above exact sequence is equivalent to:
$$
0 \to T(-2) \to (r+2)\Oc \oplus \Oc (1) \to E  \to 0
$$
A vector bundle of this type is an extension by a trivial bundle of a rank $3$
bundle defined like above, but with $r=3$.

(iii) When $\deg (Y) = 5$, as $Y$ is contained in a complete intersection $Z$ of 
type $(2,3)$, it is linked in $Z$ to a line. According to \cite{PS} 
Prop. 2.5, $Y$ has a resolution of the type:
$$
0 \to 2\Oc (-4) \to 2 \Oc (-3) \oplus \Oc(-2) \to I _Y \to 0
$$
Using this and the exact sequence \eqref{rank3} one gets that $E$ has a resolution
of the shape:
$$
0 \to 2\Oc (-1)  \xrightarrow{\varphi}  (r+1)\Oc \oplus \Oc (1) \to E \to 0
$$

Conversely, the property $\mathbf G.$ shows that, for $r \ge 3$, a general $\varphi$ 
is an injection of vector bundles, and so $E$ is a vector bundle of rank $r$. It is 
easy to show that any $E$ like above is an extension of one of rank $3$ (given by a 
sequence as above with $r=3$) by a trivial bundle.

(iv) When $\deg (Y)=6$, $Y$ is a complete intersection of type $(2,3)$, the bundle
has a resolution:
$$
0 \to \Oc (-2) \to r\Oc \oplus \Oc (1) \to E \to 0
$$ 
Again, $r \ge 3$, and any such bundle is an extension of one like above with $r=3$,
by a trivial one.

{\bf 2.  Assume $Y$ in \eqref{rank3} not connected.}

Observe that $Y$ cannot have a line $L$ as a component, because in that case 
$\Lambda _Y$ would have as summand $\Lambda _L = \Oc _L (-1)$ which has no sections.
So $\deg (Y) \ge 4$. On the other hand, when $\deg (Y) =6$, $Y$ should be 
a complete intersection, whence connected. It remains $\deg (Y) =4 {\rm \ or\ } 5$.

When $\deg (Y)=4$, $Y$ is a union of two conics. This cannot occur, because 
the minimal resolution of $Y$ has the shape:
$$
0 \to \Oc (-6) \to 2\Oc (-5) \oplus 2\Oc (-4) \to \Oc (-4) \oplus 2\Oc (-3)
\oplus \Oc (-2) \to I_Y \to 0
$$
and $I_Y(3)$ is not globally generated.

When $\deg(Y)=5$, $Y$ is a union of a conic and a cubic (plane or skew).
Again, $I_Y(3)$ is not globally generated.

(C) {\it Assume $H^0(E(-1))=0$.} Then $Y$ is contained in a complete intersection
of type $(3,3)$ and is not contained in a quadric.

{\bf 1. Assume $Y$ in \eqref{rank3} connected.} As all smooth curves of degree
$ \le 4$ are contained in a quadric, we have to analyse $5 \le \deg (Y) \le 9$.

(i) When $\deg (Y) =5$ there are three possibilities: rational quintic, elliptic
quintic, genus $2$ quintic. The last case is excluded by the fact that any genus
2 quintic lies on a quadric.

{\bf If $Y$ is a rational quintic not on a quadric}, then there is a unique 
$4$-secant $L$ to $Y$ (cf. \cite{Mi}), p. 148). According to \cite{No} there is a 
complete intersection $Z$ of type $(3,3)$ which contains $Y$ and a degree $4$ 
structure $Y'$ on $L$, which contains the first infinitesimal neighbourhood 
$L^{(2)}$. These structures are studied in \cite{BF1} (cf. \cite{BF2} and 
\cite{Ma2}), namely there is an integer $s$ and an exact sequence:
\begin{equation}\label{4str}
0 \to \Oc _L(s) \to \Oc _{Y'} \to \Oc _{L^{(2)}} \to 0 \ \ \ .
\end{equation}
The Hilbert polynomial  of $Y'$ is $\chi _{Y'}(t)=4t+s+2$. 
The liaison of $Y$ and $Y'$ in $Z$ produce the exact sequence:
\begin{equation}\label{quintic}
0 \to \omega _{Y'}\otimes \omega _Z^{-1}\ (\cong \omega _{Y'}(-2)\ ) \to \Oc _Z \to 
\Oc _Y \to 0
\end{equation}
From $ \chi_Z=\chi _Y +\chi _{\omega _{Y'}}(-2)$ follows $s=0$.

The exact sequence \eqref{quintic} produces the exact sequence:
$$
0\to I_Z(3) \to I_Y(3) \to \omega _{Y'}(1) \to 0 \ \ \ .
$$
Dualizing  \eqref{4str}, and twisting by $\Oc (1)$ one gets:
$$
0 \to \omega _{L^{(2)}}(1) \to \omega _{Y'}(1) \to \omega _L(1)\ (\cong 
\Oc _L(-1)) \to 0
$$
The last two exact sequences show that $\omega _{Y'}(1)$ and $I_Y(3)$ are not 
globally generated.

{\bf If $Y$ is an elliptic quintic} then the minimal resolution of $Y$ is of the 
shape:
$$
0 \to \Oc (-5) \to 5\Oc (-4) \to 5\Oc (-3) \to I_Y \to 0
$$
The bundles produced via \eqref{rank3} will have a resolution of the shape:
$$
0 \to \Oc(-2) \to 5\Oc(-1) \to 7\Oc \to E \to  0
$$
and consequently has Chern classes $c_1=3$, $c_2=5$, $c_3=5$. The above resolution
can be reduced to:
$$
0 \to T(-2) \oplus \Oc(-1)\to 7\Oc \to E \to 0 \ .
$$

(ii) When $\deg (Y)=6$ the only case when $I_Y$ is generated by cubics
and $Y$ is not contained in a quadric is a type of genus $3$ curves, with minimal
resolution (cf. \cite{Na}) :
$$
0 \to 3\Oc (-4) \to 4 \Oc (-3)\to I_Y \to 0 \ \ \ ,
$$
i.e. Y is arithmetically Cohen-Macaulay.
The bundles given by extensions \eqref{rank3} with such $Y$ have resolutions:
$$
0 \to 3\Oc (-1)\to 6\Oc \to E \to 0 \ \ \ ,
$$
and consequently $E$ has Chern classes $c_1=3$, $c_2=6$, $c_3=10$.

{\it Remark.} Although the case of elliptic sextics is excluded by \cite{Na},
we give also another argument:

By Serre-Hartshorne correspondence, an elliptic sextic corresponds to 
vector bundles $G$ with $c_1(G)=4$, $c_2(G)=6$, given by extensions:
\begin{equation}\label{ellsext}
0 \to \Oc \to G \to I_Y(4) \to 0 \ \ \ .
\end{equation}
$F:=G(-2)$ has $c_1(F)=0$, $c_2(F)=2$ and $H^0(F(-1))=0$, so that $F$ is stable.
These vector bundles where studied by Hartshorne (\cite{Ha}), \textsection 9. 
Namely, it is shown there that $H^0(F(1))\not =0$ and the zero set of any section
$s$ of $F(1)$ is contained in a unique smooth quadric which does depend only on $F$
and has type (3,0) on this quadric. Then $Z(s)$ has a resolution of the shape:
$$
0 \to 2\Oc (-5) \to 6\Oc(-4) \to 4\Oc(-3) \oplus \Oc (-2) \to I_{Z(s)}\to 0\ \ \ ,
$$
and, taking into account the exact sequence:
$$
0 \to \Oc \to F(1) \to I_{Z(s)}(2) \to 0 \ \ \ ,
$$
$F(1)$ has a resolution:
$$
0 \to 2\Oc (-3) \to 6\Oc(-2) \to 4\Oc(-1) \oplus 2\Oc  \to F(1)\to 0\ \ \ 
$$
Taking into account \eqref{ellsext} the elliptic sextics have resolutions:
$$
0 \to 2\Oc (-6) \to 6\Oc(-5) \to 3\Oc(-4) \oplus 2\Oc(-3)  \to I_Y\to 0\ \ \ ,
$$
(cf. also \cite{Na}) and so $I_Y(3)$ is not globally generated.

(iii) When $\deg (Y) =7$, and $Y$ is contained in a complete intersection of type 
$(3,3)$, we have to find those $Y$ linked to a degree $2$ curve in $\P^3$.

{\bf Assume $Y$ is linked to a conic (smooth, or degenerated).} Then $Y$ has a 
resolution (cf. \cite{PS}, Prop. 2.5) of the shape:
$$
0 \to \Oc (-5)\oplus \Oc (-4) \to 3\Oc (-3) \to I_Y \to 0 \ \ \ 
$$
(so that $Y$ is a curve of degree $7$ and genus $5$), 
and then the usual procedure gives vector bundles with a resolution:
$$
0\to \Oc (-2)\oplus \Oc (-1) \to 5\Oc \to E \to 0
$$
and extensions by trivial bundles of such bundles.

{\bf Remark.} It is easy to see that we obtained already all curves which are 
arithmetically Cohen-Macaulay and whose ideal are minimally generated by polynomials
of degree 3: they are the above ones, of degree $7$ and genus $5$ respectively of 
degree $6$ and genus $3$.

{\bf Assume $Y$ is linked to a disjoint union of two lines or  to a double line.}
We show in the following that this situation does not occur. 
Indeed, by \cite{GM}, Thms. 3.1, 3.2 the shape of a minimal resolution of a smooth 
curve $I_Y$ whose ideal is generated by cubics and is not arithmetically 
Cohen-Macaulay is one of the following\medskip:
$$
 0\to (a-1)\Oc(-6)\to a\Oc(-5)\oplus a\Oc(-4)\to (a+2)\Oc(-3) \to I_Y \to 0
$$
if $h^1(I_Y(2))=a-1>1$,\medskip

or
$$
 0\to (6-d)\Oc(-5)\to (15-2d)\Oc(-4)\to (10-d)\Oc (-3) \to I_Y \to 0
$$
if $h^1(I_Y(2))=0$ and $\deg(Y)=d$.\medskip

In the first case the Hibert polynomial of $Y$ is $H_Y(t)=(9-2a)t+5a-9$,
in the second it is $H_Y(t)=dt+10-2d$, and never the degree can be $7$.

(iv) When $\deg(Y)=8$, and $Y$ is contained in a complete intersection 
of type $(3,3)$, $Y$ is linked to a line. Using \cite{PS}, Prop. 2.5, $Y$ admits 
the resolution:
$$
0 \to 2\Oc(-5) \to 2\Oc(-3)\oplus \Oc (-4) \to \Oc \to \Oc_Y\to 0
$$
which produces curves of degree $8$ and genus $7$, but not with $I_Y(3)$ globally 
generated. In fact the argument from \cite{GM} used above works also in this case.

(v) Finally, if $deg(Y)=9$, $Y$ is a complete intersection. The vector bundles 
obtained in this case are extensions by trivial bundles of vector bundles given by 
exact sequences:

$$
0 \to \Oc (-3) \to 4\Oc \to E \to 0
$$

{\bf 2.  In the assumption (C) assume $Y$ not connected.}
We have seen that no connected component could be a line.
Then  $\deg\ Y\ge 4$.

When $\deg\ Y=4$, $Y$ is a union of two smooth disjoint conics. One sees directly 
that $I_Y(3)$ is not generated by global sections. 

When $\deg\ Y =5$, $Y$ is a disjoint union of a conic and a smooth cubic. The case 
of plane cubic is excluded by direct computation. We show that also is excluded 
the case $Y=Q \cup C$, where $C$ is a skew cubic and $Q$ is a  smooth conic, such that
$C \cap Q= \emptyset $. Indeed, tensoring  the exact sequence 
$$0 \to I_Q\to \Oc \to \Oc_Q \to 0$$ with $\Oc_C$
 one gets the exact sequence:
$$0 \to \Tor _1(\Oc _Q,\Oc _C) \to I_Q/I_QI_C \to \Oc _C \to 
\Oc /(I_Q+I_C) \to 0\ .$$
From here $\Tor _1(\Oc _Q, \Oc _C)= I_Q \cap I_C /I_QI_C$ and as  $Q$ and $C$ are
disjoint,  $\Tor _1(\Oc _Q, \Oc _C)=0$, so that $I_Q \cap I_C=I_QI_C$. This means
$I_Y=I_QI_C$.

Tensoring the minimal resolution of $\Oc _C$ with $I_Q$ one gets:
$$
0 \to 2I_Q(-3) \to 3I_Q(-2) \to I_Y \to 0 .
$$
Standard exact sequences of cohomology give $h^0(I_Y(3))=3$ and $h^0(I_Y(4))=13$.
Hence $I_Y(3)$ is not generated by global sections.

When $deg\ Y=6$, then $Y$ is linked in a complete intersection $Z$ of two cubic
surfaces to a locally Cohen-Macaulay curve $X$ of degree $3$. This could be:

(a) a plane cubic (smooth or degenerated)

(b) a rational skew cubic

(c) a triple structure on a line

(d) a disjoint union of a conic (degenerated or not) and a line 

(e) a disjoint union of a double line and a line

(f) a union of a conic (degenerated or not) and a skew meeting line.

(g) a union of a double line and a line meeting it.

(h) three disjoint lines\bigskip

In case (a) $Y$ is a complete intersection of type $(2,3)$, already studied.

In cases (b) and (f) $Y$ is an ACM (i.e. arithmetically Cohen-Macaulay) curve with 
minimal resolution:
$$
0 \to 3\Oc (-4) \to 4\Oc(-3) \to I_Y \to 0 \ \ \ .
$$
The vector bundles constructed via $Y$will be obtained as quotients of the shape:
$$
0 \to 3\Oc (-1) \to 6 \Oc \to E \to 0\ \ \ .
$$

As $I_Y(3)$ is globally generated, then $(I_Y/I_Y^2)(3)$ is also globally generated.
These two facts allow the assumption that one of the cubic surfaces (denote it by $S$)
containing $Y$ is smooth. Let $F:=Ax+By =0$ be the cubic form which define  $S$
containing the line of ideal $I_L=(x,y)$. Denote by $x,y,u,v$ the homogeneous 
coordinates on $\P ^3$. If we write $F$ as a polynomial with coefficients in
$\C [u,v]$, namely $F=ax+by+cx^2+dxy+ey^2+f(x,y)$, the condition that $S$ is smooth
does not allow $a,b$ to have any common zeros on $L$. So 
$I_\delta=(ax+by,x^2,xy,y^2)$ is
the ideal of the only double structure on $L$ contained in $S$ and 
$I_\tau=(ax+by+cx^2+dxy+ey^2,x^3,x^2y,xy^2,y^3)$ is the ideal of the unique 
triple structure on L contained in $S$. In fact (cf. \cite{BF2} or \cite{Ma2}) one 
has the exact sequences:
\begin{align*}
\tag{$\delta$} 0 \to \Oc_L (1)\to \Oc_\delta \to & \Oc_L \to 0 \\
\tag{$\tau$} 0 \to \Oc_L(2) \to \Oc_\tau \to & \Oc_\delta \to 0
\end{align*}
and $\omega _\delta \cong \Oc_\delta(-3)$, $\omega _\tau \cong
 \Oc_\tau (-4)$. It follows that
$\omega _\tau (1)$ is not globally generated. Then $I_Y(3)$ cannot be globally generated
either. Indeed:
in general, if $X$ is linked to $Y$ in $Z$, $Z$ being complete intersection
of type $(3,3)$, one has  the exact sequence:
$$
0 \to \omega _X \otimes \omega _Z^{-1} \to \Oc _Z \to \Oc _Y \to 0 \ \ ,
$$
which can be written also:
$$
 0 \to I_Z \to I_Y \to \omega _X \otimes \omega _Z^{-1} \to 0
$$
The condition that $I_Y (3)$ is globally generated implies that 
$\omega _X \otimes \omega _Z^{-1}(3)=\omega _X (1)$ is also globally generated.

The same argument excludes the cases when $Y$ has a line as a connected component,
i.e. (d), (e) and (h).

The only case which remains is (g). Take, as above, $S$ a smooth cubic surface 
and two concurent lines $L$ and $\ell$ on it, of ideals $I_L=(x,y)$ and $I_\ell =
(x,u)$. Then, if we start with $F$ from above, the condition that $\ell$ lies on
$S$ gives $F$ of the form: $F=ax+(b_1u+b_2v)uy+cx^2+(d_1u+d_2v)xy+euy^2 + \cdots $.
Then $I_D=(ax+\beta uy,x^2,xy,y^2)$ and $I_Y=I_D\cap I_{\ell} =
(x^2,xy,uy^2,ax+\beta uy)$. Now it is easy to compute the syzygies 
of $\Oc _Y$ and see that the Hilbert polynomial of $Y$ is
$\chi _Y(t)=3t+2$. Observe that the same result could be obtained also from the 
exact sequence:
$$
0 \to \Oc _Y \to \Oc_D \oplus \Oc_\ell \to \Oc _{D\cap \ell} \to 0 \ \ \ .
$$
$Y$ can be interpreted as a scheme which gives a locally algebraic linkage 
of $L\cup \ell$ and $L$ (cf. \cite{Ma2} or \cite{Ma3}). As such one should have 
an exact sequence of the form:
$$
0\to \Oc_L (s) \to \Oc _Y \to \Oc _{L\cup \ell} \to 0 \ \ \ .
$$ 
Taking Euler characteristics, one get $s=0$.
Dualizing the above sequence, we have a surjection $\omega _Y \to \omega _L$.
As $\omega _L (1)$ is not globally generated, $\omega _Y(1)$ is neither.

If $\deg(Y)=7$, then $Y$ is linked in a complete intersection $Z$ of type
$(3,3)$ to a locally Cohen=Macaulay curve  of degree $2$. This can be:

(a) a conic (degenerated or not)

(b) a double line

(c) two skew lines.

In the case (a) $Y$ is ACM and has resolution:
$$
0 \to \Oc(-5) \oplus \Oc(-4) \to 3\Oc(-3) \to I_Y \to 0
$$

and the corresponding rank $3$ vector bundles have resolutions of the shape:
$$
0 \to \Oc(-2)\oplus \Oc(-1) \to 5\Oc  \to E \to 0 \ .
$$
The case (b) cannot occur because for a double line on a smooth cubic surface 
$\omega _Y =\Oc _Y(-3)$ and so $\omega _Y(1)$ is not globally generated. The case 
(c) is not acceptable from the same reason (it was, in fact, already dismissed).

When $\deg(Y)=8$ then $Y$ is linked to a line and so $Y$ is ACM, with resolution of 
the shape:
$$
0 \to 2 \Oc (-5) \to \Oc(-4)\oplus 2\Oc(-3)\to I_Y \to 0
$$
and $I_Y(3)$ is not globally generated.

When $\deg(Y)=9$, then $Y$ is a complete intersection and the vector bundles 
associated to it have resolution of the type:
$$
0 \to \Oc(-3) \to 4\Oc \to E \to 0
$$

\begin{rem}
 (i) The Chern classes (Chern polynomial) of the rank $3$ vector bundles of the theorem 
(respectively the same twisted bu $-1$) proved in this paper are:
\begin{align*}
\tag{1} c(E) & = 1+ 3h+9h^2+27h^3\ , & c(E(-1)) & = 1+6h^2+20h^3 \\
\tag{2} c(E) & = 1+ 3h+6h^2+12h^3\ , & c(E(-1)) & = 1+3h^2+8h^3\\
\tag{3} c(E) & = 1+ 3h+4h^2+4h^3\ , & c(E(-1)) & = 1+h^2+2h^3\\
\tag{4} c(E) & = 1+ 3h+3h^2+3h^3\ , & c(E(-1)) & = 1+2h^3\\
\tag{5} c(E) & = 1+ 3h+7h^2+15h^3\ , & c(E(-1)) & = 1+4h^2+10h^3\\
\tag{6} c(E) & = 1+ 3h+5h^2+7h^3\ , & c(E(-1)) & = 1+2h^2+4h^3\\
\tag{7} c(E) & = 1+ 3h+6h^2+10h^3\ , & c(E(-1)) & = 1+3h^2+6h^3\\
\tag{8} c(E) & = 1+ 3h+4h^2+2h^3\ , & c(E(-1)) & = 1+h^2\\
\tag{9} c(E) & = 1+ 3h+5h^2+5h^3\ , & c(E(-1)) & = 1+2h^2+2h^3\\
\end{align*}
(ii) The Chern polynomials (1), (3), (8) decompose as follows:
\begin{align*}
\tag{1} c(E) & = (1+3h)(1+9h^2) \\
\tag{3} c(E) & = (1+2h)(1+h+2h^2)\\
\tag{8} c(E) & = (1+h)(1+2h+2h^2)\\
\end{align*}
(iii) The above decompositions correspond to extensions of vector bundles:
\begin{align*}
\tag{1} 0 \to  E' \to & E \to  \Oc(3) \to 0 \\
\tag{3} 0 \to  E' \to & E \to  \Oc(2) \to 0\\
\tag{8} 0 \to  \Oc(1) \to & E \to  TV(4) \to 0        
\end{align*}
where in each case $E'$ is a rank $2$ vector bundle (generalizations of the 
null-correlation bundle) and $TV$ is a  Trautmann-Vetter bundle (cf. \cite{T},
\cite{V}, \cite{KPR}) known also, modulo a twist due to notation,
as Tango bundle, (cf. \cite{Tango})  .
\end{rem}

\Proof (i) and (ii) are standard computations with Chern classes. 

For (iii): In the exact sequence defining $E$ in case $(1)$, let $\varphi =  
{}^t(f_1,f_2,f_3,f_4)$, where $f_1,\ldots ,f_4$ are homogeneous polynomials
of degree $3$. Let $\widetilde{\varphi}:=(-f_2,f_1,-f_4,f_3): 4\Oc \to \Oc(3)$.
Then $\widetilde{\varphi}$ is surjective and $\widetilde{\varphi}\varphi =0$, so that
$E'=\hbox{ker} \widetilde{\varphi}/\hbox{im} \varphi$ defines a generalized  
null-correlation bundle.  Standard arguments give a extension as in (iii) $(1)$.

In case (iii) $(3)$, let $\varphi = {}^t(f_1,f_2,f_3,f_4) : \Oc(-1) \to 2\Oc \oplus 
2\Oc (1)$, where $f_1,f_2$ are homogeneous polynomials of degree $1$ and $f_3,f_4$
homogeneous polynomials of degree $2$. Take $\widetilde{\varphi}:=
(-f_3, -f_4,f_1,f_2):2\Oc \oplus 2\Oc(1) \to \Oc(2)$. In this setting the proof
goes like above and one gets an exact sequence $(3)$.

(iii)$(8)$: The construction of Trautmann-Vetter (cf. \cite{T}, \cite{V}, \cite{KPR}),
applied to $\P ^3$, starts from 
the observation that $\Omega (2)$ is globally generated by $5$ global sections,
so that one gets a rank $2$ vector bundle $TV$ as a kernel:
$$
0 \to TV \to 5\Oc(-2) \to \Omega ^1  \to 0 \ \ .
$$

As $TV ^\vee \cong TV(6)$, dualizing the above exact sequence and tensoring by 
$\Oc(-2)$ one gets:
$$
0 \to T(-2) \to 5\Oc \to TV(4) \to 0 \ \ .
$$
By the defining exact sequence $h^0(E(-1))=1$ and $h^0(E(-2))=0$. 
Then, as $c_3(E(-1))=0$, a nonzero section of $E(-1)$ vanishes nowhere and one 
gets an injection of vector bundles $\Oc(1) \to E$. Let $E''$ be the cokernel of it. 
This completes to the commutative diagram with exact columns and rows:
$$
\begin{CD}
@. @. 0 @. 0 @. \\
@. @. @VVV @VVV @.\\
@. @. \Oc(1) @= \Oc(1) @. \\
@. @. @VVV @VVV @.\\
0 @>>> T(-2) @>>>5\Oc \oplus \Oc(1) @>>> E @>>> 0 \\
@. @| @VVV @VVV @.\\
0 @>>> T(-2) @>>>5\Oc  @>>> E'' @>>> 0 \\   
@. @. @VVV @VVV @.\\   
  @. @. 0 @. 0 @.      
\end{CD}
$$
where from $E'' \cong TV(4)$ .
\qed
\begin{rem}
 On $\P^3$ the Trautmann-Vetter (Tango) bundle and the null-correlation bundle 
differ only by a twist. We preffered a name or the other according to the 
construction considered. \qed
\end{rem}

\end{document}